\documentclass[12pt]{amsart}
\usepackage{amssymb}
\usepackage{times,a4wide}
\usepackage{setspace}
\newcommand{\remove}[1]{ }

\newtheorem{theorem}{Theorem}[section]
\newtheorem{proposition}[theorem]{Proposition}
\newtheorem{lemma}[theorem]{Lemma}
\newtheorem{corollary}[theorem]{Corollary}
\newtheorem{claim}[theorem]{Claim}

\theoremstyle{definition}
\newtheorem{definition}[theorem]{Definition}

\theoremstyle{remark}
\newtheorem{remark}[theorem]{Remark}

\newtheorem{examples}[theorem]{Examples}

\newcommand{\C}{\mathbb{C}}

\newcommand{\N}{\mathbb{N}}
\newcommand{\R}{\mathbb{R}}

\renewcommand{\Im}{\operatorname{Im}}
\renewcommand{\Re}{\operatorname{Re}}

\title[Geometric conditions for Interpolation]
{Geometric conditions for Interpolation in weighted spaces of entire functions}

\author[ M. Ouna\"{\i}es]
{ Myriam Ounaies}

\address{Institut de Recherche Math\'ematique Avanc\'ee, Universit\'e 
Louis Pasteur 7 Rue Ren\'e Des\-car\-tes, 67084 Strasbourg CEDEX, France.}
\email{ounaies@math.u-strasbg.fr}

\date{\today}



\begin{document}

\begin{abstract} 
We use $L^2$ estimates for the $\bar\partial$ equation to find geometric conditions on discrete interpolating varieties for  weighted spaces $A_p(\C)$ of entire functions such that $\vert f(z)\vert\le Ae^{Bp(z)}$ for some $A,B>0$.  In particular, we give a characterization when $p(z)=e^{\vert z\vert}$ and more generally when $\ln p(e^r)$ is convex and $\ln p(r)$ is concave.
\end{abstract}

\maketitle

\section*{Introduction}

Let $p$ be a weight (see Definition \ref{defweight} below) and  $A_p(\C)$ be the vector space of  entire functions satisfying $\displaystyle \sup_{z\in \C}\vert f(z)\vert e^{-Bp( z)}<\infty$ for some $B>0$.  

For  instance, if $p(z)=\vert z\vert$ then $A_p(\C)$ is the space of all entire functions of exponential growth. More generally when $p(z)=\vert z\vert^\alpha$, $\alpha>0$ $A_p(\C)$ is the space of entire functions of order $<\alpha$ or of order $\alpha$ and of finite type.
When $p(z)=\log(1+\vert z\vert^2)+\vert \Im z\vert$, $A_p(\C)$ is the space of Fourier transforms of distributions with compact support in the real line. 

We are concerned with the interpolation problem for $A_p(\C)$. That is, finding conditions on a given discrete sequence of complex numbers $V=\{z_j\}_j$ so that, for any  sequence of complex numbers $\{w_j\}_j$ with convenient growth conditions, there exists $f\in A_p(\C)$ such that $f(z_j)=w_j$, for all $j$. We will then say that $V$ is an interpolating variety for the weight $p$.
We actually consider the problem with prescribed multiplicities on each $z_j$, but for the sake of simplicity, we will assume the multiplicities to be equal to $1$ in the introduction.

There exists an analytic characterization of interpolating varieties for all weights $p$ satisfying Definition \ref{defweight} (see \cite{Be-Li1}).

We are interested in finding a geometric description which would enable us to decide whether a discrete sequence is interpolating for $A_p(\C)$ by looking at the density of the points.
This was done for the weight $p(z)=\log(1+\vert z\vert^2)+\vert \Im z\vert$ in \cite{Ma-Or-Ou}. In the present work we will mainly treat radial weights.

The geometric conditions will be given in terms of $N(z,r)$, the integrated counting function of the points of $V$ in the disk of center $z$ and radius $r$ (see Definition \ref{count} below).

When $p$ is radial $(p(z)=p(\vert z\vert))$ and doubling $(p(2z)\le 2 p(z))$, Berenstein and Li \cite{Be-Li} gave a geometric characterization of interpolating varieties for $p$, namely,
\begin{itemize}
\item[(i)]  $N(z_j,\vert z_j\vert)=O( p(z_j))$ when $j\rightarrow \infty$;
\item [(ii)] $N(0,r)=O(p(r))$ when $r\rightarrow \infty$.
\end{itemize}

Hartmann and Massaneda (\cite[Theorem 4.3]{Ha-Ma}) gave a proof of this theorem based on $L^2$ estimates for the solution to the $\bar\partial$-equation of the sufficiency provided that $p(z)=O( \vert z\vert^2\Delta p(z))$. Note that we can always regularize $p$ into a smooth function, see Remark  \ref{reg} below.

In this paper we will give a proof in the same spirit as \cite{Ha-Ma} without the assumption on the Laplacian of $p$ (see Theorem \ref{doub}).

When the above condition  on the Laplacian is satisfied, we will prove that (i) is necessary and sufficient (see Theorem \ref{finiteorder}).

In \cite{Be-Li}, Berenstein and Li also studied rapidly growing radial weights, allowing infinite order functions in $A_p(\C)$, as $p(z)=e^{\vert z\vert}$, and more generally weights such that $\ln p(e^r)$
is convex. They gave sufficient conditions as well as necessary ones.

We will give a characterization of interpolating varieties for the weight $p(z)=e^{\vert z\vert}$ and more generally for weights $p$ such that 
$p(z)=O(\Delta p(z))$ (see Theorem \ref{exp}) and also for radial $p$ when $\ln p(e^r)$ is convex and $\ln p(r)$ is concave for large $r$ (see Theorem \ref{convex}).

In particular, we will show that $V$ is interpolating for $A_p(\C)$, $p(z)= e^{\vert z\vert}$, if and only if 
$$N(z_j,e)=O(e^{\vert z_j\vert}), \hbox{ when } \ j\rightarrow \infty.$$

The difficult part  in each case is the sufficiency. As in \cite{Be-Or, Ha-Ma, Ma-Or-Ou}, we will follow a Bombieri-H\"ormander approach based on $L^2$-estimates on the solution to the $\bar\partial$-equation. The scheme will be the following: the condition on the density gives a smooth interpolating function $F$ with a good growth such that the support of $\bar\partial F$ is  far from the points $\{z_j\}$ (see Lemma \ref{F}). Then we are led to solve the $\bar\partial$-equation: $\bar\partial u=-\bar\partial F$ with $L^2$-estimates, using a theorem of H\"ormander \cite{Ho1}. To do so, we need to construct a subharmonic function $U$ with a convenient growth and with prescribed singularities on the points $z_j$. Following Bombieri \cite{Bo}, the fact that $e^{-U}$ is not summable near the points  $\{z_j\}$ forces  $u$ to vanish on the points $z_j$ and we are done by defining the interpolating entire function by $u+F$.

The delicate point of the proof is the construction of the function $U$. It is done in two steps: first we construct a function $U_0$ behaving like $ \ln \vert z-z_j\vert^2$ near $z_j$ with a good growth and  with a control on $\Delta U_0$ (the Laplacian of $U_0$), thanks to the conditions on the density and the hypothesis on the weight itself. Then we add a function $W$ such that $\Delta W$ is large enough so that $U=U_0+W$ is subharmonic.

A final remark about notation:

$A$, $B$ and $C$ will denote positive constants and 
their actual value may change from one  occurrence to the next. 

$F(t)\lesssim G(t)$ means that there exists constants $A,B>0$, not depending on $t$ such that $F(t)\le A G(t)+B$ while 
$F\simeq G$ means that $F\lesssim G\lesssim F$.

The notation $D(z,r)$ will be used for  the Euclidean disk of center $z$ and radius $r$. We will write $\displaystyle \partial f=\frac{\partial f}{\partial z}$, $\displaystyle \bar \partial f=\frac{\partial f}{\partial \bar z}$. Then  $\Delta f=4\partial\bar\partial f$ denotes the Laplacian of $f$. 

To conclude the introduction, the author wishes to thank X. Massaneda for useful talks and remarks.

\section{Preliminaries and main results}

\begin{definition}\label{defweight}
A subharmonic function $p:\C\longrightarrow\R_+$,  is called a weight 
if
\begin{itemize}
\item[(a)] $\ln(1+|z|^2)=O(p(z))$;
\item[(b)]  there exist constants $C_1, C_2>0$ such that $\vert z-w\vert \le 1$ implies $p(w)\le C_1p(z)+C_2$.
\end{itemize}
\end{definition}

Note that condition (b) implies that $p(z)=O(\exp (A\vert z\vert)$ for some $A>0$.

We will say that the weight is ``radial'' when $p(z)=p(\vert z\vert)$ and that it is ``doubling'' when $p(2r)\lesssim  p(r)$. 

Let $A(\C)$ be the set of all entire functions, we consider the space
$$A_p(\C)=\Bigl\{f\in A(\C):\ \  \forall z\in \C,\  |f(z)| \le A\, e^{B p(z)}\hbox{ for some } A>0, B>0\Bigr\}.$$

\begin{remark}\label{rem}{}\mbox{}

\begin{itemize}
\item[(i)] Condition (a) implies that $A_p(\C)$ contains all polynomials.
\item[(ii)] Condition (b) implies that $A_p(\C)$ is stable under differentiation.
\end{itemize}
\end{remark}

\begin{examples}\mbox{}
\begin{itemize}
\item $p(z)=\ln(1+\vert z\vert^2)+\vert \Im z\vert$. Then $A_p(\C)$ is the space of Fourier transforms of distributions with compact support in the real line.
\item $p(z)=\ln(1+\vert z\vert^2)$. Then $A_p(\C)$ is the space of all the polynomials.
\item $p(z)=\vert z\vert$. Then $A_p(\C)$ is  the space of entire functions of exponential type.
\item $p(z)=\vert z\vert^{\alpha}$, $\alpha>0$. Then $A_p(\C)$ is the space of all entire functions of order $< \alpha$ or of order $\alpha$ and finite type.
\item $p(z)=e^{\vert z\vert^\alpha}$, $0<\alpha\le 1$.
\end{itemize}
\end{examples}

Throuought the paper $V=\{(z_j,m_j)\}_{j\in \N}$ will denote a multiplicity variety, that is, a sequence of points $\{z_j\}_{j\in \N}\subset \C$ such that $\vert z_j\vert \rightarrow \infty$, and a sequence of positive integers $\{m_j\}_{j\in \N}$ corresponding to the multiplicities of the points $z_k$ and $p$ will denote a weight. 

\begin{definition}
We  say that $V$ is an interpolating variety for $A_p(\C)$ if for every doubly indexed sequence $\{w_{j,l}\}_{j,0\le l<m_j}$ of complex numbers such that, for some positive constants $A$ and $B$ and for all $j\in \N$, 
 $$\sum_{l=0}^{m_j-1}\vert w_{j,l}\vert\le A e^{Bp(z_j)},$$ we can find an entire function $f\in A_p(\C)$, 
with
$$\frac{f^{l}(z_j)}{l!}=w_{j,l}$$
for all $j\in \N$ and $0\le l<m_j$.
\end{definition}

\begin{remark}(See \cite[Proposition 2.2.2]{Be-Ga}) 
Thanks to condition (b), we have, for some constants $A>0$ and $B>0$,  
 $$\forall z\in \C,\ \ \sum_k \left\vert\frac{f^{(k)}(z)}{k!}\right\vert\le A e^{Bp(z)}.$$ 
  \end{remark}

If we consider the space
$$A_p(V)=\Bigl\{W=\{w_{j,l}\}_{j,0\le l < m_j}\subset \C: \ \  \forall j,\ \sum_{l=0}^{m_j-1}|w_{j,l}| \le A\, e^{B p(z_j )} \hbox{ for some } A>0, B>0\Bigr\}$$
and we define the restriction map by
\[ 
\begin{split}
\mathcal R_V: A_p & (\C)\longrightarrow A_p (V)\\
& f\quad \mapsto\; \left\lbrace \frac{f^{l}(z_j)}{l!}\right\rbrace _{j,0\le l\le m_j-1}, 
\end{split} 
\]
we may equivalently define the interpolating varieties by the property that  $\mathcal R_V$ maps $A_p(\C)$ onto $A_p(V)$.

Note that $A_p(\C)$ can be seen as the union of the Banach spaces 
$$A_{p,B}(\C)=\{f\in A(\C), \ \ \Vert f\Vert_B=\sup_{z\in \C}\vert f(z)\vert e^{-Bp(z)} <\infty\}$$ 
and has a structure of an (LF)-space
with the the inductive limit topology.
The same can be said about $A_p(V)$.

The problem we are considering is to find conditions on $V$ so that it is an interpolating variety for $A_p(\C)$. 

In order to state the geometric conditions, we define the counting function and the integrated counting function:

\begin{definition}\label{count}
For $z\in \C$ and $r>0$ we set
\[
n(z,r)=\sum\limits_{|z-z_j|\le r} m_j
\]
and
\[
\begin{split}
N(z,r) & =\int_0^r\frac{n(z,t)-n(z,0)}t\, dt + n(z,0)\ln r \\
& =\sum_{0<\vert z- z_j\vert\le r}m_j \ln \frac{r}{\vert z-z_j\vert}+n(z,0)\ln r .
\end{split}
\]
\end{definition}

\begin{remark}\label{reg}
The weight $p$ may be regularized as in \cite[Remark 2.3]{Ha-Ma}  by replacing $p$ by its average over the disc $D(z,1)$.
 Thus we may suppose $p$ to be of class ${\mathcal C}^2$ when needed.
 \end{remark}

Before stating our results we recall that if $V$ is interpolating for $A_p(\C)$, then the  following two conditions are necessary:
 \begin{equation}\label{N(z,1)}
\exists A>0,\  \exists B>0,\ \ \forall j\in \N,\ \ \ N(z_j,e) \le A\ p(z_j)+B
\end{equation}
and when $p$ is radial, 
\begin{equation}\label{N(r)}
\exists A>0,\  \exists B>0,\ \ \forall r>0,\ \ N(0,r)\le A\,p(r)+B.
\end{equation}
See Theorem \ref{gen} and Proposition \ref{gen1} for the proof.

\smallskip 
In condition (i), we may replace $N(z_j,e)$ by $N(z_j, c)$ with any constant $c>1$ .
We are now ready to state our main results. 
We begin by giving  sufficient conditions for a discrete variety $V$ to be interpolating for radial weights.

\begin{theorem}\label{suff}
Assume the weight $p$ to be radial.
If condition \eqref{N(z,1)} holds  and 
\begin{equation}\label{n(0,r)}
\exists A>0,\  \exists B>0\ \  \forall r>0,
\ \ \int_0^r n(0,t)dt\le A\ p(r)+B
\end{equation}
 then $V$ is interpolating for $A_p(\C)$.
\end{theorem}
We  note that 
$$\int_0^r n(0,t)dt= \sum_{\vert z_j\vert\le r} m_j(r-\vert z_j\vert )\le rN(0,r).$$

Consequently, by condition (2), we see that $\displaystyle \int_0^r n(0,t)dt \le A r p(r)+B$, for some $A>0$ and $B>0$ is a necessary condition.

Adapting  our method to the doubling case we find the characterization given by Berenstein and Li \cite[Corollary 4.8]{Be-Li}:

\begin{theorem}\label{doub}
Assume $p$ to be radial and  doubling. 

\smallskip

$V$  is interpolating for $A_p(\C)$ if and only if conditions \eqref{N(r)} holds and   

\begin{equation}\label{N(z,z)}
\exists A>0,\  \exists B>0\ \ \forall j\in \N,\ \ \ N(z_j,\vert z_j\vert) \le A\ p(z_j)+B.
\end{equation}

\end{theorem}

The theorem holds if we replace $N(z_j,\vert z_j\vert)$ by $N(z_j,C \vert z_j\vert)$ for any constant $C>0$.
Note that radial and doubling weights satisfy $p(r)=O(r^\alpha)$ for some $\alpha>0$. In other words, they have at most a polynomial growth.
Examples of radial and doubling weights are $p(z)=\vert z\vert^\alpha (\ln (1+\vert z\vert^2))^\beta$, $\alpha>0, \beta\ge 0$, but for $p(z)=\vert z\vert^\alpha$, we can give a better result:

\begin{theorem}\label{finiteorder}
Assume that  $p(z)=O( \vert z\vert^2\Delta p(z))$ and \par

(b') $\exists C_1>0$, $\exists C_2>0$,\ such that $\vert z-w\vert \le \vert z\vert$ implies $p(w)\le C_1p(z)+C_2$.

\smallskip 

$V$ is interpolating for $A_p(\C)$ if and only if condition \eqref{N(z,z)} holds.
 
\end{theorem}

\begin{remark}\label{rad}
It is easy to see that radial and doubling weights  satisfy condition (b').

Theorem  \ref{finiteorder} applies to $p(z)=\vert z\vert^\alpha$, $\alpha>0$. 
For this weight and  with the extra assumption that  there is a function $f\in A_p(\C)$ vanishing on every $z_j$ with multiplicity $m_j$, it was shown in (\cite[Theorem 3]{Sq1}) that condition
 \eqref{N(z,z)} is sufficient and necessary. 
 \end{remark}
  \smallskip

Next we are interested in the case where $p$ grows rapidly, allowing infinite order functions in $A_p(\C)$. A fundamental example is $p(z)=e^{\vert z\vert}$.

In \cite{Be-Li}, Berenstein and Li studied this weight and more generally those for which $\ln p(e^r)$ is convex. 
They gave sufficient conditions (Corollaries 5.6 and Corollary 5.7)  as well as necessary ones (Theorem 5.14, Corollary 5.15).  

The following result gives a characterization in particular for the weight $p(z)=e^{\vert z\vert}$ and more generally for $p(z)=e^{\vert z\vert^\alpha}, \ \alpha\ge 1$.

\begin{theorem}\label{exp}
Assume that $p(z)=O(\Delta p(z))$.

\smallskip
 $V$ is interpolating for $A_p(\C)$ if and only if condition \eqref{N(z,1)} holds.
\end{theorem}

The next theorem will give a characterization when $p$ is radial, $q(r)=\ln p(e^r)$ is convex and  $\frac{r}{q'(\ln r)}=\frac{p(r)}{p'(r)}$ is increasing (for large $r$). If we set
$u(r)=\ln p(r)$, we have $\frac{p(r)}{p'(r)}=\frac{1}{u'(r)}$. Thus, the last condition means that $u(r)$ is concave for large $r$.
 We recall that  the convexity of $q$  implies that $p(r)\ge A r+B$, for some $A,B>0$ (see \cite[Lemma 5.8]{Be-Li}).

The weights $p(z)=\vert z\vert^\alpha$, $\alpha>0$ and $p(z)=e^{\vert z\vert}$ satisfy these conditions. Examples of weights also satisfying these conditions but not those of Theorems \ref{doub}, \ref{finiteorder} or \ref{exp} are $p(z)=e^{\vert z\vert^\alpha}$, $0<\alpha<1$ and $p(z)=e^{[\log(1+\vert z\vert^2)]^\beta}$, $\beta>1$.

\begin{theorem}\label{convex}
Assume that $p$ is a radial weight and that for a certain $r_0>0$ it satisfies the following properties
\begin{itemize}
\item $q(r):=\ln p(e^r)$ is convex on $[\ln r_0,\infty[$;
\item $q'(\ln r_0)>0$ and $\frac{r}{q'(\ln r)}$ is increasing on $[r_0, \infty[$.
\end{itemize}
Then $V$ is interpolating for $A_p(\C)$ if and only if the following condition holds:
 \begin{equation}\label{conv}
\exists A>0, \exists B>0, \ \ \forall \vert z_j\vert\ge r_0, N(z_j,\max(\frac{\vert z_j\vert}{q'(\ln \vert z_j\vert)},e))\le Ap(z_j)+B.
\end{equation}
\end{theorem}

The theorem holds if we replace $\displaystyle \frac{\vert z_j\vert}{q'(\ln \vert z_j\vert)}$ by $\displaystyle \frac{C\vert z_j\vert}{q'(\ln \vert z_j\vert)}$ for any constant $C>0$.  

When $p(z)=\vert z\vert^\alpha$, conditions \eqref{conv} and \eqref{N(z,z)} are the same and when $p(z)=e^{\vert z\vert}$, conditions \eqref{conv} and \eqref{N(z,1)} are the same.

As immediate corollaries of Theorem \ref{convex}, we have the following

\begin{corollary}
Let  $p(z)=e^{\vert z\vert^\alpha}$, $0<\alpha< 1$.
$V$ is interpolating for $A_p(\C)$ if and only if the following condition holds:
 \begin{equation}
\exists A>0, \exists B>0, \ \ \forall j, \ \ N(z_j,\vert z_j\vert^{1-\alpha})\le A p(z_j)+ B.
\end{equation}
\end{corollary}

\begin{corollary}
Let  $p(z)=e^{[\ln(1+\vert z\vert^2)]^\beta}$, $\beta\ge 1$.
$V$ is interpolating for $A_p(\C)$ if and only if the following condition holds:
 \begin{equation}
\exists A>0, \exists B>0, \ \ \forall j, \ \ N(z_j,\vert z_j\vert [\ln(1+\vert z_j\vert^2)]^{1-\beta})\le A p(z_j)+ B.
\end{equation}
\end{corollary}

\section{General results about the interpolation theory}

\smallskip 

For the sake of completeness, we include in this section some standard results about interpolating varieties. 
Let us begin by well known necessary conditions.

\begin{theorem}\label{gen}\cite[Theorem 1]{Sq1}
Assume that  $V$ is an interpolating variety for $A_p(\C)$. Let $R_j$ be positive numbers satisfying
\begin{equation}\label{Rj}
\vert z-z_j\vert\le R_j \Longrightarrow p(z)\le C_1 p(z_j)+C_2
\end{equation}
where $C_1$ and $C_2$ are positive constants not depending on $j$. Then  the following condition holds: 
\begin{equation}\label{N(z,r)}
\exists A>0,\ \exists B>0,\ \forall j,\ \ N(z_j,R_j)\le A\ p(z_j)+B.
\end{equation}
\end{theorem}
\begin{remark}\label{exam}
In view of property (b) of the weight $p$, $R_j=1$ satisfies condition \eqref{Rj}. 
In the case where $p$ satisfies the stronger condition (b') (as in Theorem \ref{finiteorder}) we are allowed the larger numbers $R_j=\vert z_j\vert$. 
\end{remark} 

\begin{proposition}\label{gen1}
Assume that  $V$ is an interpolating variety for $A_p(\C)$. Then the following condition is satisfied
\begin{equation}\label{N(0,r)}
\exists A>0,\ \exists B>0,\ \forall R>0,\ \ N(0,R)\le A\ \max_{\vert z\vert =R} p(z)+B.
\end{equation}
\end{proposition}
The proof of this proposition is implicitely contained in \cite[page 26]{Be-Li}.

\begin{definition}
We say that $V$ is weakly separated if there exist constants $A>0$ and $B>0$ such that
\begin{equation}\label{sep}
\frac{1}{\delta_j^{m_j}}\le A e^{B p( z_j) }
\end{equation}
for all $j\in \N$, where 
 $$\delta_j:= \inf \left\lbrace 1,\frac{1}{2}\inf _{k\not=j}|z_j-z_k|\bigr)\right\rbrace $$
 is called the separation radius.
\end{definition}

\begin{lemma}\label{weaksep}
If \eqref{N(z,1)} holds  then $V$ is weakly separated.
\end{lemma}

\begin{proof}
Fix $j\in \N$ and let $z_l\not=z_j$ be such that $\vert z_j-z_l\vert =\inf_{k\not= j} |z_j-z_k|$. If $\vert z_j-z_l\vert\ge 2$, then $\delta_j=1$. Otherwise, $2\delta_j=\vert z_j-z_l\vert$ and the following inequalities hold
$$N(z_l,e)\ge \sum_{0<\vert z_k-z_l\vert\le e}m_k\ln\frac{e}{\vert z_k-z_l\vert}\ge m_j\ln\frac{e}{\vert z_j-z_l\vert}=m_j\ln\frac{e}{2\delta_j}\ge \ln\frac{1}{\delta_j^{m_j}}.$$
Since by condition \eqref{N(z,1)} and property (b) of the weight,
$$N(z_l,e)\le A +Bp(z_l)\le A +Bp(z_j),$$
we readily deduce the desired estimate.
\end{proof}
	
We will follow the same scheme as in \cite{Ha-Ma, Ma-Or-Ou}, first constructing a smooth interpolating function with the right growth.

\begin{lemma}\label{F}
Suppose $V$ is weakly separated. 
Given $W=\{w_{j,l}\}_{j\in \N,0\le l \le m_j-1} \in  A_p(V)$, there exists a smooth function $F$ such that
\begin{itemize}
\item for some constants $A>0$ and $B>0$, $\vert F(z)\vert \le Ae^{Bp(z)}$, $\vert \bar\partial F(z)\vert \le Ae^{Bp(z)}$ for all $z\in \C$; 
\item the support of $\bar\partial F$
 is contained in the union of the annuli
\[
A_j=\{z\in\C:\  \frac{\delta_j}{2}  \le|z-z_j|\le
 \delta_j\} ;
\] 
\item $\displaystyle \frac{F^{(l)}(z_j)}{l!}=w_{j,l}$\  for all $j\in \N$, $0\le l\le m_j-1$.
\end{itemize}
\end{lemma}

A suitable function $F$ is of the form 
\[
F(z)=\sum_j  \mathcal X\Bigl(\frac{|z-z_j|^2}
{\delta_j^2}\Bigr)\sum\limits_{l=0}^{m_j-1} w_{j,l} (z-z_j)^l,
\]
where $\mathcal X$ is a smooth cut-off function with  $\mathcal X(x)=1$ if $|x|\le1/4$ and $\mathcal X(x)=0$ if $|x|\ge 1$. See \cite{Ha-Ma} or \cite{Ma-Or-Ou} for details of the proof.

Now, when looking for a holomorphic interpolating function of the form
$f=F+u$, we are led to the $\bar\partial$-problem
\[
\bar\partial u=-\bar\partial F\ , 
\]
which we solve using H\"ormander's theorem \cite[Theorem 4.2.1]{Ho2}.

The interpolation problem is then reduced to the following: 

\begin{lemma}\label{U}

If $V$ is weakly separated and if there exists a subharmonic function $U$ satisfying for  certain constants $A,B>0$,

\begin{itemize}
\item[(i)] $U(z)\le Ap(z)+B$ for all  $z\in \C $;
\item[(ii)] $-U(z)\le Ap(z)+B$ for $z$ in the support of $\bar\partial F$;
\item[(iii)] $U(z)\simeq m_j\ln |z-z_j|^2$ near $z_j$,
\end{itemize}
then $V$ is interpolating.
\end{lemma}

\begin{proof}[Proof.]\par
The weak separation gives a smooth interpolating function $F$ (see Lemma \ref{F}). From H\"ormander's theorem \cite[Theorem 4.4.2]{Ho1}, we can find a ${\mathcal{C}}^\infty$ function $u$ such that 
$\bar\partial u=-\bar\partial F$ and, denoting by $d\lambda$ the Lebesgue measure,
\begin{equation}\label{hor}
\int_{\C}\frac{\vert u(w) \vert^2e^{-U(w)-Ap(w)}}{
(1+\vert w\vert^2)^2}\,d\lambda(w)\le\int_{\C}\vert \bar\partial F \vert^2e^{-U(w)-Ap(w)}\,d\lambda(w).
\end{equation}
By  property (a) of the weight $p$, there exists $C>0$ such that
$$\int_{\C}e^{-Cp(w)}d\lambda(w)<\infty.$$
Thus, using (ii) of the lemma, and the estimate on $|\bar\partial F(z)|^2$, we see that  the last integral in \eqref{hor} is convergent if  $A$ is large enough.
By condition (iii), near $z_j$, $e^{-U(w)}(w-z_j)^l$ is not summable for $0\le l\le m_j-1$, so we have necessarily $u^{(l)} (z_j)=0$ for all $j$ and $0\le l \le m_j-1$ and consequently, 
$\displaystyle \frac{f^{(l)}(z_j)}{ l!}=w_{j,l}$.

Now, we have to verify that $f$ has the desired growth. By the mean value inequality,
$$
\vert f(z)\vert\lesssim \int_{D(z,1)}\vert f(w)\vert\,d\lambda(w)\lesssim
\int_{D(z,1)}\vert F(w)\vert\,d\lambda(w)
+\int_{D(z,1)}\vert u(w)\vert\,d\lambda(w).
$$
Let us estimate the two integrals that we denote by $I_1$ and $I_2$. For $w\in D(z,1)$, 
$$\vert F(w)\vert\lesssim e^{Bp( w)}\lesssim e^{Cp( z)}.$$
Hence 
$$I_1\lesssim e^{C p(z)}.$$
\smallskip

To estimate $I_2$ we use Cauchy-Schwarz inequality
$$I_2^2\le J_1\,J_2$$
where
$$
J_1=\int_{D(z,1)}\vert u(w)\vert^2e^{-U(w)-Bp( w)}\,d\lambda(w),\ J_2=\int_{D(z,1)}e^{U(w)+Bp(w)}\,d\lambda(w).$$
By property (a) of $p$ we have 
$$J_1 \lesssim \int_{\C} \vert u(w) \vert^2e^{-U(w)-Bp(w)}\,d\lambda(w) \lesssim \int_{\C}\frac{\vert u(w) \vert^2e^{-U(w)-Bp(w)}}{
(1+\vert w\vert^2)^2}\,d\lambda(w)< +\infty$$
if $B>0$ is chosen big enough.

To estimate $J_2$ we use  condition (i) of the lemma and  property (b) of the weight $p$. For $w\in D(z,1)$,  
$$e^{U(w)+Bp(w)}\le e^{Cp( w)}\lesssim e^{Ap(z)}.$$
We easily deduce that $J_2\lesssim e^{Ap(z)}$ and finally that $f\in A_p(\C)$.
\end{proof}

\section{Proofs of the main theorems}

We will use a smooth cut-off function $\mathcal X$ with $\mathcal X(x)=1$ if $\vert x\vert \le 1/4$ and $\mathcal X(x)=0$ if $\vert x\vert \ge 1$.

\begin{remark}\label{finite}
In the proofs of the sufficiency part, we may need to assume that for all $j$, we have $\vert z_j\vert\ge a$ for a suitable $a>0$. This will be done without loss of generality up to a linear transform and in view of property (b) of the weight.
\end{remark}

\begin{proof}[Proof of Theorem \ref{suff}]${}$\par

By Lemma \ref{weaksep}, condition \eqref{N(z,1)} implies the weak separation. So we are done if we construct a function $U$ satisfying the conditions of Lemma \ref{U}.

 Set $\mathcal X_j(z)=\mathcal X(\vert z-z_j\vert^2)$.

In order to construct the desired function we begin by defining 
$$U_0(z)=\sum_j m_j \mathcal X_j (z) \ln{\vert z-z_j\vert^2}.$$
Note that there is locally a finite number of non vanishing terms in the sum and that each term (and consequently $U_0$) is nonpositive.
It is also clear that $U_0(z)- m_j\ln\vert z-z_j\vert^2$ is continuous near $z_j$.

 We want to estimate $-U_0$ on the support of $\bar \partial F$, and the ``lack of subharmonicity'' of $U_0$, then we will add a correcting term to obtain the function $U$ of the lemma.

Suppose $z$ is in the support of $\bar \partial F$. We want to show that $-U_0(z)\lesssim p(z)$.
Let  $k$ be the unique integer such that $\frac{\delta_k}{ 2 } \le |z-z_k| \le \delta_k$. Then

$$-U_0(z)\le 2\sum_{\vert z-z_j\vert\le 1}  m_j  \ln \frac{1 }{ \vert z-z_j\vert}=2m_k \ln\frac{1}{ \vert z-z_k\vert}+2\sum_{j\not=k,\vert z-z_j\vert\le 1} m_j  \ln \frac{1 }{ \vert z-z_j\vert}.$$

Using that $\vert z-z_k\vert\ge \frac{\delta_k}{ 2}$ and that for $j\not=k$ we have

$$|z_k-z_j| \le |z-z_j| + |z-z_k| \le 2 |z-z_j|,$$
we obtain that
\begin{equation}\label{ii}
-U_0(z)\le  2\ln\frac{1}{ \delta_k^{m_k}}+2 N(z_k,2)\lesssim p(z_k)\lesssim p(z).
\end{equation}
The last inequalities follows from condition \eqref{N(z,1)}, the weak separation \eqref{sep} and property (b) of the weight $p$.

Now we want to get a lower bound on $\Delta U_0(z)$. We have 
\[
\begin{split}
\Delta U_0(z) & =\sum_j m_j \mathcal X_j(z)\Delta \ln \vert z-z_j\vert^2 \\
& +8 \Re \left(\sum_j m_j \bar\partial \mathcal X_j(z) \partial \ln \vert z-z_j\vert^2\right)+4 \sum_j m_j \partial \bar\partial \mathcal X_j (z)  \ln \vert z-z_j\vert^2.
\end{split}
\]
The first sum is a positive measure and on the supports of $\bar\partial  \mathcal X_j$ and $ \partial \bar\partial \mathcal X_j$, we see that $1/2\le \vert z-z_j\vert\le 1$. Consequently, for a certain constant $\gamma>0$ we have

\begin{equation*}\label{Delta U_0}
\Delta U_0(z)\ge - \gamma (n(z,1)-n(z,1/2))\ge -\gamma n(z,1) \ge -\gamma (n(0,\vert z\vert+1)-n(0,\vert z\vert-1)).
\end{equation*}

We set $n(0,t)=0$ if $t<0$, 
$$ f(t)=\int_{t-1}^{t+1} n(0,s)ds,\  \ g(t)=\int_0^t f(s) ds \hbox{ and }\  W(z)=g(|z|).$$

Let us compute the Laplacian of $W$, taking the derivatives in the sense of distributions
$$\Delta W(z)=\frac{1}{ \vert z\vert} g'(\vert z\vert)+g''(\vert z\vert )\ge g''(\vert z\vert)=f'(\vert z\vert )=n(0,\vert z\vert+1)-n(0,\vert z\vert-1).$$
The function $U$ defined by
$$U(z)=U_0(z)+\gamma W(z)$$
is then clearly subharmonic. 
On the other hand, using condition \eqref{n(0,r)} and (b) we have the following inequalities: 
\[
f(s)\le 2 n(0,s+1),\ \ W( z)=g(\vert z\vert)\le 2 \int_1^{\vert z\vert+1} n(0,s) ds \lesssim  p(z).
\]
Since $U_0\le 0$ it is clear that $U$ satisfies condition (i) of Lemma \ref{U}.
Using the estimate \eqref{ii} and the fact that $W$ is nonnegative we see that $U$ satisfies  condition (ii). Finally as condition (iii) is also already fulfilled by $U_0$ and  $W$ is continous, it is also fulfilled by $U$. 
\end{proof}

\begin{proof}[Proof of Theorem \ref{doub}]${}$\par

{\it Necessity.}
In view of Remark \ref{rad} and Remark \ref{exam}, we apply Theorem \ref{gen} with $R_j=\vert z_j\vert$ and we readily obtain the necessity of  \eqref{N(z,z)}.

Condition \eqref{N(r)} is necessary  by  Proposition \ref{gen1}.

{\it Sufficiency.}
By Lemma \ref{weaksep}, condition \eqref{N(z,z)} implies the weak separation.  We will proceed as in Theorem \ref{suff}, constructing a function $U$ satisfying (i), (ii) and (iii) from  Lemma \ref{U}. Thanks to the doubling condition, we can control the weight $p$ in discs $D(z_j,\vert z_j\vert)$ instead of just $D(z_j,e)$ in the general case. We will construct $U_0$ as in the previous theorem, except that we now take $\mathcal X_j$'s with supports of radius $\simeq \vert z_j\vert$:

Set $$\mathcal X_j(z)=\mathcal X\left(\frac{16\vert z-z_j\vert^2}{ \vert z_j\vert^2}\right),$$
and introduce the negative function
$$U_0(z)=\sum_j m_j \mathcal X_j (z) \ln\frac{16\vert z-z_j\vert^2}{ \vert z_j\vert^2}.$$

When $z$ is in the support of $\bar \partial F$, let  $k$ be the unique integer such that $\frac{\delta_k}{ 2 } \le |z-z_k| \le \delta_k$. 
Repeating the estimate on $-U_0(z)$, we have

$$-U_0(z)\le 2\sum_{\vert z-z_j\vert\le \frac{\vert z_j\vert}{ 4}}  m_j  \ln \frac{\vert z_j\vert }{ 4\vert z-z_j\vert}\le 2m_k \ln\frac{\vert z_k\vert }{ \delta_k}+2\sum_{0<\vert z_k-z_j\vert\le \frac{\vert z_j\vert}{ 2}} m_j  \ln \frac{\vert z_j\vert }{ 2\vert z_k-z_j\vert}.$$
We have $\displaystyle \frac{\vert z_j\vert}{ 2}\le \vert z_k\vert$ whenever $\displaystyle \vert z_k-z_j\vert \le \frac{\vert z_j\vert }{ 2}$. We deduce the inequalities
\begin{equation}\label{-U0}
-U_0(z)\le 2\ln \frac{1}{ \delta_k^{m_k}}+ 2N(z_k,\vert z_k\vert)\lesssim p(z_k)\lesssim p(z).
\end{equation}
Again, the last inequalities follow from condition \eqref{N(z,1)}, the weak separation \eqref{sep} and property (b) of the weight $p$.

We estimate  $\Delta U_0(z)$ as before except that now $\vert \bar\partial \mathcal X_j(z)\vert \lesssim \frac{1}{ \vert z_j\vert}$ and $\vert \partial\bar\partial \mathcal X_j(z)\vert\lesssim \frac{1}{ \vert z_j\vert^2}$. On the support of these derivatives, $\displaystyle \frac{\vert z_j\vert }{ 8}\le \vert z-z_j\vert\le \frac{\vert z_j\vert}{ 4}$ 
and $\displaystyle \frac{\vert z\vert}{ 2}\le \vert z_j\vert\le 2\vert z\vert$. We deduce that  
\[
\Delta U_0(z)\gtrsim - \frac{n(0,2\vert z\vert)-n(0,\frac{\vert z\vert}{ 2})}{ \vert z\vert^2}.
\] 
To construct the correcting term $W$  set
$$ f(t)=\int_0^t n(0,s)ds,\  \ g(t)=\int_0^t \frac{f(s)}{ s^2} ds \hbox{ and }\  W(z)=g(2 |z|).$$
Finally,  to estimate the Laplacian of $W$, we write $t=2|z|$. We have
$$\Delta W(z)=4\left( \frac{1}{ t} g'(t)+g''(t)\right) =\frac{4}{ t^2} \left( f'(t)-\frac{f(t)}{ t}\right) $$
and
$$f(t)=\int_0^t n(0,s)ds = \int_0^{\frac{t}{4}} n(0,s)ds+\int_{\frac{t}{4}}^t n(0,s)ds\le\frac{t}{ 4}n\left( 0,\frac{t}{4}\right) +t\left( 1- \frac{1}{4}\right)  n(0,t).$$
Thus, $$f'(t)-\frac{f(t)}{ t}=n(0,t)-\frac{f(t)}{ t} \ge\frac{1}{ 4}\left( n(0,t)-n\left( 0,\frac{t}{4}\right) \right) $$
and $$\Delta W(z)\gtrsim \frac{ n(0,2|z|)-n(0,\frac{|z|}{ 2})}{ |z|^2}.$$
Now, the desired function will be of the form 
$$U(z)=U_0(z)+\gamma W(z),$$
where $\gamma $ is a positive constant sufficiently large.
The following inequalities are easy to see:
$$f(t)\le t n(0,t),\ \ g(t)\le \int_0^t \frac{n(0,s) }{ s} ds=N(0,s).$$
Thus,  by condition \eqref{N(r)} and the doubling condition, 
$$0\le W(z)\le N(0,2|z|)\lesssim p(2z)\lesssim p(z).$$
We conclude that $U$ satisfies all the desired conditions.
\end{proof}

\begin{proof}[Proof of Theorem \ref{finiteorder}]\mbox{}

{\it Necessity.} 
Recalling Remark \ref{exam}, we apply once again Theorem \ref{gen} to deduce the necessity of condition \eqref{N(z,z)}.

{\it Sufficiency.}
The proof is the same as for Theorem \ref{doub}, we only change the estimate on $\Delta U_0$ and the correcting  term $W$.
Let us have a new look at $\Delta U_0(z)$.
\begin{equation}\label{Delta U0}
\Delta U_0(z)\gtrsim -\sum_{\vert z-z_j\vert\le \frac{\vert z_j\vert}{ 4} }\frac{
m_j }{ \vert z\vert^2}.
\end{equation}
If the sum is not empty, let $z_k$ be the point appearing in the sum with the largest modulus. For all $z_j$ such that $\vert z-z_j\vert\le \frac{\vert z_j\vert}{4}$, we have
$$\vert z_j-z_k\vert \le \vert z-z_k\vert +\vert z-z_j\vert \le \frac{\vert z_k\vert}{ 4}+\frac{\vert z_j\vert }{ 4}\le \frac{\vert z_k\vert}{ 2}.$$
We deduce that
$$\Delta U_0(z)\gtrsim -\frac{ n(z_k, \frac{\vert z_k\vert}{ 2})}{ \vert z\vert^2}.$$
Besides, 
$$n(z_k,\frac{\vert z_k\vert}{ 2})\le m_k+\frac{1}{ \ln 2} \sum_{0<\vert z_j-z_k\vert\le \frac{\vert z_k\vert}{ 2}}m_j\ln\frac{\vert z_k\vert }{ \vert z_j-z_k\vert} \lesssim N(z_k,\vert z_k\vert)\lesssim p(z_k).$$
Note that $\vert z-z_k\vert \le \frac{\vert z_k\vert}{4}$ implies that $\vert z-z_k\vert \le \vert z\vert$. Thus by condition (b') we have $p(z_k)\lesssim p(z)$.
Finally we get $$\Delta U_0(z)\gtrsim -\frac{p(z)}{ \vert z\vert^2}\gtrsim -\Delta p(z).$$ 
Then we take
$$U(z)=U_0(z)+\gamma p(z)$$
where $\gamma$ is a positive constant chosen large enough.
\end{proof}

\begin{proof}[Proof of Theorem \ref{exp}]

We already know by Theorem \ref{gen} that condition \eqref{N(z,1)} is necessary.

Let us consider the  function $U_0$ that we constructed in the proof of Theorem \ref{suff}. Again, we only change the estimate on $\Delta U_0$ and the correcting term $W$. 
We find 
\begin{equation}
\Delta U_0(z)\gtrsim - n(z,1).
\end{equation}
If $n(z,1)\not=0$, let $z_k$ be in $D(z_k,1)$.
Then 
$$n(z,1)\le n(z_k,2)\le m_k+\frac{1}{ 1-\ln 2}\sum_{0<\vert z_k-z_j\vert < 2}m_j\ln \frac{e}{ \vert z_k-z_j\vert}\lesssim N(z_k,e)\lesssim p(z_k)\lesssim p(z).$$

The function 
$$U(z)=U_0(z)+\gamma p(z),$$
with $\gamma>0$ large enough has the desired properties.
\end{proof}

\begin{proof}[Proof of Theorem \ref{convex}]${}$\par

We set $c:=\inf(q'(\ln r_0),1)$ and $\displaystyle \psi(r)=\frac{r}{q'(\ln r)}$. 

\begin{claim}\label{prel}$ $ 
\begin{itemize}
\item [(i)] Let $r\ge 2r_0$. Then $c\psi(r)\le r$ and 
\[
\vert x\vert\le c \frac{\psi(r)}{ 2 } \textrm{implies that } \frac{\psi(r)}{ 2}\le \psi(r+x)\le 2 \psi(r);
\]
\item [(ii)]  For all $r\ge r_0$,
\[
p(r+\psi(r))\le e p(r).
\]
\end{itemize}
\end{claim}

Assuming this claim true for the moment, let us proceed with the proof of Theorem \ref{convex}.
{\it Necessity.}
In view of Theorem \ref{gen}, it suffices to show that $R_j=\psi(\vert z_j\vert)$ satisfy condition \eqref{Rj}. Let $\vert z_j\vert\ge r_0$ and $w$ be such that $\vert w-z_j\vert \le R_j$. Thus $\vert w\vert \le \vert z_j\vert + \psi(\vert z_j\vert)$ and as a consequence of (ii) of the claim, we obtain 
\[
p(w)\le p(\vert z_j\vert + \psi(\vert z_j\vert)\le e p(\vert z_j\vert).
\] 

{\it Sufficiency.}
We may assume that $\vert z_j\vert\ge r_0$ for all $j$ (see Remark \ref{finite}).
We apply Lemma \ref{weaksep} to deduce that $V$ is weakly separated. We repeat the proof of Theorems \ref{finiteorder} and \ref{doub}, replacing $\vert z_j\vert$ by $c\psi(\vert z_j\vert)$. More precisely, we set 
 $$\mathcal X_j(z)=\mathcal X\left(\frac{16\vert z-z_j\vert^2}{ c^2\psi(\vert z_j\vert )^2}\right)$$
and we define the negative function
$$U_0(z)=\sum_j m_j \mathcal X_j (z) \ln\frac{16\vert z-z_j\vert^2}{ c^2\psi(\vert z_j\vert)^2}.$$
We use (i) of Claim \ref{prel} to obtain that whenever $\vert z-z_j\vert \le \frac{c\psi(\vert z_j\vert)}{2}$, we have
\[
\begin{split}
\psi(\vert z\vert) &\le \psi(\vert z_j\vert+\vert z-z_j\vert)\le 2 \psi(\vert z_j\vert);\\
\frac{\psi(\vert z_j\vert)}{2} &\le \psi(\vert z_j\vert-\vert z-z_j\vert)\le \psi(\vert z\vert).
\end{split}
\]
Adapting the inequality \eqref{-U0} when $\frac{\delta_k}{2}\le \vert z-z_k\vert\le \delta_k$ and applying condition \eqref{conv} we find  
\[
-U_0(z)\le 2\ln \frac{1}{ \delta_k^{m_k}}+ 2N(z_k,c\psi(\vert z_k\vert)\lesssim p(z_k)\lesssim p(z).
\]
Let us now find a lower bound for $\Delta U_0$ on $\C$. By analogy to \eqref{Delta U0} and the inequalities following we obtain 
\[ 
\Delta U_0(z)\gtrsim -\sum_{\vert z-z_j\vert \le \psi(\vert z_j\vert)c/4}\frac{1}{\psi (\vert z_j\vert)^2}\gtrsim -\frac{N(z_l,c\psi(\vert z_l\vert))}{(\psi (\vert z\vert))^2}\,
\]
where $z_l$ is one of the points  appearing in the sum with the largest modulus. 
Recall that  
\[
\vert z-z_l\vert \le  \psi(\vert z_l\vert)c/4\le \psi (\vert z\vert)c/2\le \psi (\vert z\vert) 
\]
and consequently from (ii) of Claim \ref{prel} we deduce that 
\[
p(z_l)\le p(\vert z\vert +\vert z-z_l\vert)\le p(\vert z\vert+\psi(\vert z\vert))\le e p(z).
\]
Finally, we apply condition \eqref{conv} 
\[
\Delta U_0(z)\gtrsim -\frac{p(z)}{ \psi (\vert z\vert)^2}.
\]
Let us compute the Laplacian of $p(z)=e^{q(\ln \vert z\vert)}$ in terms of the convex function $q$. Setting $r=\vert z\vert$,
$$\Delta p(z)=\frac{p'(r)}{ r}+p''(r)=\frac{[q'(\ln r)]^2}{ r^2}p(r)+\frac{q''(\ln r)}{ r^2}p(r)\ge \frac{[q'(\ln r)]^2}{ r^2}p(r)=\frac{p(r)}{ [\psi(r)]^2}.$$
We readily deduce that
\[
\Delta U_0(z)\ge -\gamma \Delta p(z)
\]
for some $\gamma>0$.
As in the preceding proofs, the function $U(z)=U_0(z)+\gamma p(z)$ satisfies the properties stated in Lemma \ref{U}.
\end{proof} 

\begin{proof} [Proof of Claim \ref{prel}]

Let $r>r_0$. A computation gives
\[
\psi'(r)=\frac{q'(\ln(r))-q''(\ln r) }{ [q'(\ln r)]^2}.
\]
Recall that $q''$ is nonnegative. Besides as $q'$ is increasing we have $q'(\ln r)\ge q'(\ln r_0)\ge c$. We deduce that $0\le \psi'(r)\le \frac{1}{ c}$. Note that we also have the inequality $c\psi(r)\le r$.

Assume now $r\ge 2r_0$ and let $\vert x\vert \le \frac{c\psi(r)}{2}$. Then
\[
\vert r+ x\vert \ge r-\frac{c\psi(r)}{ 2}\ge \frac{r}{2}\ge r_0.
\]
We can now use the mean value theorem and the preceding estimates to write
\[
\vert \psi(r+x)-\psi(r)\vert \le \frac{\vert x\vert }{ c}\le \frac{\psi(r)}{2}.
\]
We easily deduce (i).
\smallskip

To prove (ii), put $u(t)=q(\ln t)=\ln p(t)$. We have $u'(t)=\frac{1}{\psi (t)}\le \frac{1}{ \psi(r)}$ for all $t\ge r\ge r_0$ because $\psi$ is increasing on $[r_0, \infty[$.
Applying the mean value theorem again, we obtain $u(r+\psi(r))-u(r)\le 1$. We deduce that
\[
p(r+\psi(r))=e^{u(r+\psi(r))-u(r)} p(r)\le e\, p(r).\qedhere
\]
\end{proof}

\bibliographystyle{plain}

\begin{thebibliography}{10}

\bibitem{Be-Ga}
C.~A. Berenstein and R.~Gay.
{\em Complex analysis and special topics in harmonic analysis}.
Springer-Verlag, New York, 1995.

\bibitem{Be-Li}
C.~A. Berenstein and B.~Q. Li.
\newblock Interpolating varieties for spaces of meromorphic functions.
\newblock {\em J. Geom. Anal.}, 5(1):1--48, 1995.

\bibitem{Be-Li1}
Carlos~A. Berenstein and Bao~Qin Li.
\newblock Interpolating varieties for weighted spaces of entire functions in
  {${\bf C}\sp n$}.
\newblock {\em Publ. Mat.}, 38(1):157--173, 1994.

\bibitem{Be-Or}
B.~Berndtsson and J.~Ortega~Cerd{\`a}.
\newblock On interpolation and sampling in {H}ilbert spaces of analytic
  functions.
\newblock {\em J. Reine Angew. Math.}, 464:109--128, 1995.

\bibitem{Bo}
E.~Bombieri.
\newblock Algebraic values of meromorphic maps.
\newblock {\em Invent. Math.}, 10:267--287, 1970.


\bibitem{Ha-Ma}
A.~Hartmann and X.~Massaneda.
\newblock On interpolating varieties for weighted spaces of entire functions.
\newblock {\em J. Geom. Anal.}, 10(4):683--696, 2000.

\bibitem{Ho1}
Lars H{\"o}rmander.
\newblock {\em An introduction to complex analysis in several variables}.
\newblock D. Van Nostrand Co., Inc., Princeton, N.J.-Toronto, Ont.-London,
  1966.

\bibitem{Ho2}
Lars H{\"o}rmander.
\newblock Generators for some rings of analytic functions.
\newblock {\em Bull. Amer. Math. Soc.}, 73:943--949, 1967.

\bibitem{Ma-Or-Ou}
X.~Massaneda, J.~Ortega, and M.~Ouna\"{\i}es.
\newblock A geometric characterization of interpolation in
  {$\hat{\mathcal{E}}\sp{\prime}(\mathbb{R})$}.
\newblock {\em Trans. Amer. Math. Soc.},  358(8): 3459-3472, 2006 .


\bibitem{Sq1}
W.~A. Squires.
\newblock Geometric condition for universal interpolation in {$\hat{\mathcal
  {E}}\sp{\prime} $}.
\newblock {\em Trans. Amer. Math. Soc.}, 280(1):401--413, 1983.

\end{thebibliography}

\end{document}